# SUBEXPONENTIAL ASYMPTOTICS OF HYBRID FLUID AND RUIN MODELS

By Bert Zwart, Sem Borst and Krzystof Dębicki

*Eindhoven University of Technology, CWI and University of Wrocław*

We investigate the tail asymptotics of the supremum of $X(t) + Y(t) - ct$, where $X = \{X(t), t \geq 0\}$ and $Y = \{Y(t), t \geq 0\}$ are two independent stochastic processes. We assume that the process $Y$ has subexponential characteristics and that the process $X$ is more regular in a certain sense than $Y$. A key issue examined in earlier studies is under what conditions the process $X$ contributes to large values of the supremum only through its average behavior. The present paper studies various scenarios where the latter is *not* the case, and the process $X$ shows some form of "atypical" behavior as well. In particular, we consider a fluid model fed by a Gaussian process $X$ and an (integrated) On-Off process $Y$. We show that, depending on the model parameters, the Gaussian process may contribute to the tail asymptotics by its moderate deviations, large deviations, or oscillatory behavior.

**1. Introduction.** Consider two independent stochastic processes $X = \{X(t), t \geq 0\}$ and $Y = \{Y(t), t \geq 0\}$. In the present paper we investigate the tail asymptotics of the supremum of the superposition of $X$ and $Y$ with an additional linear drift term, that is, we are interested in the behavior of

$$(1.1) \qquad \mathbb{P}\Big\{\sup_{t \geq 0}[X(t) + Y(t) - ct] > u\Big\}, \qquad u \to \infty.$$

The latter probability may be interpreted as an overflow probability in queueing theory, but also as a ruin probability. Motivated by applications in both queueing and ruin problems, we are especially interested in the case where at least one of the processes $X$ and $Y$ has *subexponential* characteristics.

A key problem which has received a lot of attention is under what conditions the process $X$ in (1.1) can be replaced by its mean, that is, under what conditions

$$(1.2) \quad \mathbb{P}\left\{\sup_{t\geq 0}[X(t) + Y(t) - ct] > u\right\} \sim \mathbb{P}\left\{\sup_{t\geq 0}[Y(t) - ct] > u\right\}, \qquad u \to \infty,$$

when $\mathbb{E}\{X(t)\} = 0$. Such an asymptotic equivalence is commonly termed a "reduced-load" equivalence. Results of the form (1.2) have been shown to hold under various assumptions; see, for instance, [[2, 8, 18]–[21, 23, 24]].

The goal of the present paper is to examine various scenarios that may arise when (1.2) does *not* hold, and the process $X$ contributes to the asymptotics in (1.1) by behaving in some "atypical" manner as well. Throughout the paper, $Y$ is assumed to be an (integrated) On-Off process: during On-periods, $Y(t)$ increases at rate $r$, while $Y(t)$ remains constant during Off-periods. The "noise" process $X$ is assumed to be Gaussian with zero mean; precise definitions and assumptions concerning $X$ and $Y$ are given in Section 2. The value of the rate $r$, relative to the negative drift $c$, has a crucial effect on the qualitative behavior of (1.1). All references listed above consider sufficient conditions for (1.1) in the case $r > c$. In the opposite case $r \leq c$, the right-hand side of (1.2) is 0 since $Y(t) \leq rt$, so that (1.2) cannot hold.

After a more detailed model description in Section 2, we will state the main result of the paper in Section 3 which considers the case $r > c$. In this case, (1.2) may or may not hold, depending on the tail behavior of the On-period $T_{\text{on}}$. If $X$ is a Brownian motion, then (1.2) holds if $\mathbb{P}\{T_{\text{on}} > u\} \sim \mathbb{P}\{T_{\text{on}} > u - \sqrt{u}\}$ as $u \to \infty$; see Theorem 3.1 in [19]. Here, we investigate what happens when this is *not* the case, and show that $X$ may then contribute to the tail behavior of $V_{X+Y}^c := \sup_{t\geq 0}[X(t) + Y(t) - ct]$ in a quite complicated fashion. Informally speaking, $X$ contributes to large values of $V_{X+Y}^c$ by its *moderate* deviations. The first part of Section 3 considers the strongly related problem of sampling a Brownian motion at a subexponential time. This part relies on recent work of Foss and Korshunov [14]. In Sections 4 and 5 we turn to the case $r \leq c$, which, as mentioned above, implies that the right-hand side of (1.2) is 0. In this case, the typical way for $V_{X+Y}^c$ to reach a large value is fundamentally different, depending on whether $r < c$ or $r = c$ (obviously, both cases differ from the case $r > c$). In the case $r = c$, which is treated in Section 4, a large value of $V_{X+Y}^c$ is associated with a single long On-period of $Y$ and oscillatory behavior of $X$. In the case $r < c$, which is studied in Section 5, the tail behavior of $V_{X+Y}^c$ is determined by the large-deviations behavior of $X$, which needs to occur during a single long On-period generated by the process $Y$.

**2. Model description and preliminaries.** We consider the supremum $V_{X+Y}^c = \sup_{t\geq 0}[X(t) + Y(t) - ct]$ of the superposition of two independent stochastic processes $X = \{X(t), t \geq 0\}$ and $Y = \{Y(t), t \geq 0\}$ with an additional drift term. We assume that $\mathbb{E}\{X(1) + Y(1)\} < c$ to ensure that $V_{X+Y}^c$ is finite a.s.

Note that $V_{X+Y}^c$ may be interpreted as the "buffer content" in steady state in a fluid queue with "drain rate" $c$, with $X(t)$ and $Y(t)$ denoting the amount of traffic generated by the processes $X$ and $Y$ during the time interval $[-t, 0]$. We will frequently make comparisons with the buffer content for each of the two processes in isolation. For $c > \mathbb{E}\{X(1)\}$ and $D \subseteq [0, \infty)$, define $V_X^c(D) := \sup_{t \in D}[X(t) - ct]$, and let $V_X^c := V_X^c([0, \infty))$ be a random variable representing the steady-state workload in a buffer with drain rate $c$ fed by the process $X$ only; $V_Y^c(D)$ and $V_Y^c$ are defined similarly.

Before specifying the processes $X$ and $Y$ in more detail, we first introduce some useful notational conventions and concepts. If $T$ is a nonnegative random variable with finite mean, then $T^r$ is a random variable with density $\mathbb{P}\{T > x\}/\mathbb{E}\{T\}$. For any two real functions $f(\cdot)$ and $g(\cdot)$, we use $f(u) \sim g(u)$ to denote that $f(u) = g(u)(1 + o(1))$ as $u \to \infty$, that is, $\lim_{u \to \infty} f(u)/g(u) = 1$. We further write $f(u) \lesssim g(u)$ to indicate that $\limsup_{u \to \infty} f(u)/g(u) \leq 1$. We use various classes of distributions. In particular, we consider the class $\mathcal{L}$ of long-tailed distributions, the class $\mathcal{S}$ of subexponential distributions and the class $\mathcal{R}$ of regularly varying distributions. We also consider the subclass $\mathcal{S}^*$ of $\mathcal{S}$. For definitions and further background on these classes, we refer to [13].

Throughout the paper, the process $X(t)$ is a centered [i.e., $\mathbb{E}\{X(t)\} = 0$] Gaussian process with stationary increments, a.s. continuous sample paths, $X(0) = 0$ a.s., and variance function $\mathrm{Var}\{X(t)\} = \sigma_X^2(t)$. We impose the following conditions in Sections 4 and 5:

C1. $\sigma_X^2(t) \in C([0, \infty))$ is increasing;
C2. $\sigma_X^2(t)$ is regularly varying at 0 with index $\beta \in (0, 2]$ and $\sigma_X^2(t)$ is regularly varying at $\infty$ with index $\alpha \in (0, 2)$.

Two important examples satisfying C1 and C2 are (i) fractional Brownian motion with Hurst parameter $H \in (0, 1)$ [in this paper denoted by $B_H(t), t \geq 0$] and (ii) the class of *integrated Gaussian processes* which has been studied extensively in [10].

The process $Y(t)$ is an (integrated) On-Off process with stationary increments. For future use we give an explicit construction of such a process, following [16]: Let $\{T_{\mathrm{on},m}, m \geq 0\}$ be a sequence of i.i.d. random variables representing the On-periods of the source. Similarly, let $\{T_{\mathrm{off},m}, m \geq 0\}$ be the Off-periods; generic On- and Off-periods are denoted by $T_{\mathrm{on}}$ and $T_{\mathrm{off}}$. Define further an independent random variable $I$ such that

$$p = \mathbb{P}\{I = 1\} = \frac{\mathbb{E}\{T_{\mathrm{on}}\}}{\mathbb{E}\{T_{\mathrm{on}}\} + \mathbb{E}\{T_{\mathrm{off}}\}} = 1 - \mathbb{P}\{I = 0\}.$$

To obtain a stationary alternating renewal process, we define the delay random variable $D_0$ by

$$D_0 = I T_{\mathrm{on},0}^r + (1 - I)(T_{\mathrm{off},0}^r + T_{\mathrm{on},0}).$$

Then the delayed renewal sequence

$$\{Z_n, n \geq 0\} = \left\{ D_0, D_0 + \sum_{m=1}^{n} (T_{\mathrm{off},m} + T_{\mathrm{on},m}), n \geq 1 \right\}$$

is stationary.

Next, we define the process $\{J(t), t \geq 0\}$ as follows. $J(t)$ is the indicator of the event that the source is On at time $t$. Formally, we have

$$J(t) = I \mathbb{1}_{\{t < T_{\mathrm{on},0}^r\}} + (1 - I) \mathbb{1}_{\{T_{\mathrm{off},0}^r \leq t < T_{\mathrm{off},0}^r + T_{\mathrm{on},0}\}} + \sum_{n=0}^{\infty} \mathbb{1}_{\{Z_n + T_{\mathrm{off},n+1} \leq t < Z_{n+1}\}}.$$

The On-Off process $\{J(t), t \geq 0\}$ is strictly stationary, see Theorem 2.1 of [16]. The process $\{Y(t), t \geq 0\}$ is then defined by

$$Y(t) := r \int_0^t J(s) \, ds.$$

Note that the mean rate of $Y(t)$ is $\rho = pr$. We assume $\rho < c$ to ensure that $V_{X+Y}^c$ is finite a.s.

For the distribution of the On-period $T_{\mathrm{on}}$ we impose the following condition in Section 3 (a similar condition has been introduced by Borovkov [6] in a related problem):

T1. The tail of the random variable $T$ has the form $\mathbb{P}\{T > u\} = e^{-L(u)u^\beta}$, with $0 < \beta < 1$, and $L(u)$ slowly varying and twice differentiable. Moreover, $L'(t) = o(L(t)/t)$ and $L''(t) = o(L(t)/t^2)$.

The next result shows that a random variable $T$ is subexponential if it satisfies T1. In fact, one can derive a slightly stronger result:

LEMMA 2.1. *If $T$ satisfies* T1, *then $T \in \mathcal{S}^*$. In particular, $T, T^r \in \mathcal{S}$.*

PROOF. The hazard function $Q$ and hazard rate $q$ of $T$ are given by $Q(u) = L(u)u^\beta$ and

$$q(u) = \beta L(u) u^{\beta - 1} + u^\beta L'(u).$$

Hence, we have $q(u) \to 0$, $uq(u) \to \infty$ and $uq(u)/Q(u) \to \beta \in (0, 1)$. According to Corollary 3.9 of [15], this implies that $T \in \mathcal{S}^*$, which in turn implies $T, T^r \qquad\qquad\qquad\qquad \in \qquad\qquad\qquad\qquad \mathcal{S}.$ $\qquad\square$

**3. Moderately heavy tails and moderate deviations.** In this section we assume that $Y(t), t \geq 0$, is an integrated On-Off process with peak rate $r > c$ and that $X(t) = W(t), t \geq 0$, is a standard Brownian motion. Under this condition, following Theorem 3.2 in [19], the reduced-load equivalence (1.2) then holds if

$$(3.1) \qquad \mathbb{P}\{T_{\text{on}} > u\} \sim \mathbb{P}\{T_{\text{on}} > u - \sqrt{u}\}.$$

If $T_{\text{on}}$ has a Weibullian tail of the form $e^{-u^\beta}$, then (3.1) and (1.2) both hold if $\beta < 1/2$. Moreover, neither (3.1) nor (1.2) holds if $\beta > 1/2$. If (3.1) does not hold, then we call $T_{\text{on}}$ *moderately heavy-tailed*, following [3]. This section aims to obtain the tail asymptotics of $V^c_{X+Y}$ in the moderately heavy-tailed case.

Define $W_\mu(t) = W(t) + \mu t, t \geq 0$. Our main result is the following:

THEOREM 3.1. *If $X(t) = W(t)$ is a standard Brownian motion, and $Y(t)$ is an integrated On-Off process with peak rate $r > c$, with $T_{\text{on}}$ satisfying* T1, *then*

$$\mathbb{P}\{V^c_{X+Y} > u\} \sim p\frac{r-\rho}{c-\rho}\mathbb{P}\{W_{r-c}(T^r_{\text{on}}) > u\}.$$

To prove Theorem 3.1, we can exploit the fact that $W$ has independent increments. This naturally leads to the framework of Asmussen, Schmidli and Schmidt [4]: the increment process of $S(t) := W(t) + Y(t) - ct$ is regenerative, with regeneration points being the ends of On-periods. Thus, the analysis consists of two steps, which are carried out in Sections 3.1 and 3.2. In Section 3.1 we investigate the tail behavior of $W_{r-c}(T)$, with $T$ an independent subexponential random variable. We believe that the results in this section are of independent interest. After that, we apply the results of the first step to obtain the tail behavior of $V^c_{X+Y}$ in Section 3.2, culminating in a proof of Theorem 3.1.

3.1. *Sampling a Brownian motion at a subexponential time.* Let $T$ be a random variable which is long-tailed and independent of $\{W_\mu(t), t \geq 0\}$. Define the running maximum $M_\mu(t) = \sup_{0<s<t} W_\mu(s)$. The goal of this section is to determine the tail behavior of $W_\mu(T)$ for $\mu > 0$. As the first step towards that goal, we show that $W_\mu(T)$ and $M_\mu(T)$ are tail-equivalent. The following lemma establishes this tail equivalence under minimal assumptions.

LEMMA 3.1. *If $T \in \mathcal{L}$, then $W_\mu(T), M_\mu(T) \in \mathcal{L}$, and $\mathbb{P}\{W_\mu(T) > u\} \sim \mathbb{P}\{M_\mu(T) > u\}$.*

PROOF.   Define for $x > 0$, $\tau(x) := \inf\{t : W_\mu(t) = x\}$ and fix $y > 0$. Note that $\tau(x+y) \overset{d}{=} \bar{\tau}(x) + \bar{\tau}(y)$, with the latter two random variables distributed as $\tau(x)$ and $\tau(y)$, but mutually independent. Write, for some $M$ and $K$,

$$\frac{\mathbb{P}\{M_\mu(T) > x + y\}}{\mathbb{P}\{M_\mu(T) > x\}} = \frac{\mathbb{P}\{T > \bar{\tau}(x) + \bar{\tau}(y)\}}{\mathbb{P}\{T > \tau(x)\}}$$

$$\geq \mathbb{P}\{\tau(y) < M\} \frac{\mathbb{P}\{T > \tau(x) + M\}}{\mathbb{P}\{T > \tau(x)\}}$$

$$\geq \mathbb{P}\{\tau(y) < M\} \int_K^\infty \frac{\mathbb{P}\{T > z + M\}}{\mathbb{P}\{T > z\}} \, d\mathbb{P}\{\tau(x) \leq z\}.$$

Note that, since $T \in \mathcal{L}$, we can choose (for each $\varepsilon > 0, M < \infty$) an appropriate constant $K = K(M, \varepsilon)$ such that

$$\frac{\mathbb{P}\{T > z + M\}}{\mathbb{P}\{T > z\}} \geq 1 - \varepsilon$$

when $z \geq K$. Combining these observations gives the lower bound

(3.2)        $$\frac{\mathbb{P}\{M_\mu(T) > x + y\}}{\mathbb{P}\{M_\mu(T) > x\}} \geq (1 - \varepsilon)\mathbb{P}\{\tau(y) < M\}\mathbb{P}\{\tau(x) > K\}.$$

Noting that $\tau(x) \to \infty$ a.s. as $x \to \infty$, we obtain the property $M_\mu(T) \in \mathcal{L}$ by letting first $x \to \infty$ in (3.2) and then $\varepsilon \to 0$, $M \to \infty$.

Next, observe that

$$\mathbb{P}\{M_\mu(T) > x\}$$

$$\leq \mathbb{P}\{W_\mu(T) \geq x - y\} + \mathbb{P}\{W_\mu(T) < x - y; M_\mu(T) > x\}$$

$$\leq \mathbb{P}\{W_\mu(T) \geq x - y\} + \mathbb{P}\{\tau(x) < T; W_\mu(T) - W_\mu(\tau(x)) < -y\}$$

$$\leq \mathbb{P}\{W_\mu(T) \geq x - y\} + \mathbb{P}\left\{\tau(x) < T; \inf_{t > \tau(x)}[W_\mu(t) - W_\mu(\tau(x))] < -y\right\}$$

$$= \mathbb{P}\{W_\mu(T) > x - y\} + \mathbb{P}\{M_\mu(T) > x\}\mathbb{P}\{V_{W_0}^\mu > y\},$$

where the last equality follows from the strong Markov property for $W_\mu(t)$. We conclude that

$$\mathbb{P}\{W_\mu(T) > x - y\} \geq \mathbb{P}\{M_\mu(T) > x\}\mathbb{P}\{V_{W_0}^\mu \leq y\}.$$

From this inequality, the obvious property

$$\mathbb{P}\{W_\mu(T) > x - y\} \leq \mathbb{P}\{M_\mu(T) > x - y\},$$

and the fact that $M_\mu(T) \in \mathcal{L}$, one obtains the tail equivalence of $W_\mu(T)$ and $M_\mu(T)$, and, in particular, the property $W_\mu(T) \in \mathcal{L}$.   □

We now examine the tail behavior of $\mathbb{P}\{W_\mu(T) > u\}$ as $u \to \infty$ in the moderately heavy-tailed regime. A related problem has been investigated by Foss and Korshunov [14]: they consider the random variable $N(T)$, with $N(\cdot)$ a renewal process. As their analysis shows, the computations in the moderately heavy-tailed regime are very technical. We could apply a similar approach here (using explicit formulas for Brownian motion and the Laplace method), but we will follow a different approach: we construct a renewal process $N_\mu(t)$ with the property

$$(3.3) \qquad M_\mu(t) - 1 \leq N_\mu(t) \leq M_\mu(t),$$

which, in view of Lemma 3.1, reduces the problem to the one studied in [14]. This approach avoids a lot of tedious computations and may be of independent interest.

We construct the renewal process $N_\mu(t)$ as follows. Define a sequence of stopping times $\tau_i$, $i \geq 1$, by

$$\tau_i := \inf\{t : W_\mu(t) = i\}.$$

Then, define

$$N_\mu(t) := \max\{n : \tau_n \leq t\}.$$

It is obvious that (3.3) holds. Moreover, $N_\mu(t)$ is a renewal process, since $\tau_i - \tau_{i-1}$, $i \geq 1$, is an i.i.d. sequence. Define

$$\Lambda(x) = \sup_y [xy - \log \mathbb{E}\{e^{y\tau_1}\}].$$

Let $\lambda(x)$ be the optimizing point in the above supremum. Since

$$\mathbb{E}\{e^{y\tau_1}\} = e^{\mu - \sqrt{\mu^2 - 2y}},$$

we have

$$\Lambda(x) = \frac{\mu^2}{2}x - \mu + \frac{1}{2x} \quad \text{and} \quad \lambda(x) = \frac{\mu^2}{2} - \frac{1}{2x^2}.$$

We now state the main result of this section.

PROPOSITION 3.1. *If $T$ satisfies* T1, *then*

$$\mathbb{P}\{W_\mu(T) > u\} \sim \mathbb{P}\{M_\mu(T) > u\}$$
$$\sim \mathbb{P}\{N_\mu(T) > u\} \sim e^{-H(t(u), u)},$$

*with $H(t, u) = Q(t) + u\Lambda(t/u)$ and $t(u)$ a solution of $Q'(t) = -\lambda(t/u)$.*

PROOF.   Assumption T1 implies that $Q(u) = -\log \mathbb{P}\{T > u\}$ is twice differentiable and that $uQ''(u) \to 0$. This allows us to apply Theorem 5.1 of [14] to obtain the tail behavior of $N_\mu(T)$. The remaining assertions follow from Lemma 3.1, (3.3) and the fact that $T \in \mathcal{L}$.   □

If $T$ has a Weibullian tail, that is, $Q(u) = u^\beta$, $0 < \beta < 1$, then Lemma 6.3 of [14] implies

$$(3.4) \qquad u/\mu - t(u) \sim uQ'(u/\mu) = \beta u^\beta \mu^{1-\beta}.$$

This indicates that a large value of $W_\mu(T)$ is caused by a realization of $T$ which is about $u/\mu - \beta u^\beta \mu^{1-\beta}$, which implies that $W_\mu(u) - \mu u$ must be of the order $u^\beta$.

Hence, if $1/2 < \beta < 1$ [in which case the asymptotic equivalence $\mathbb{P}\{W_\mu(T) > u\} \sim \mathbb{P}\{\mu T > u\}$ does not hold], then $W_\mu(u) - \mu u$ contributes to the asymptotics by means of its *moderate deviations*. If $\beta < 1/2$, then the deviation of $W_\mu(u)$ from its mean falls within the fluctuations at the CLT level, in which case $\mathbb{P}\{W_\mu(T) > u\} \sim \mathbb{P}\{\mu T > u\}$. This can be shown using Proposition 3.1 (see [14]) or, directly, by invoking results from [19].

We conclude this section with two results which are of crucial importance in the next section. The first result concerns the question whether or not the tail distribution of $W_\mu(T)$ is subexponential.

PROPOSITION 3.2.   *If $T$ satisfies* T1, *then $W_\mu(T) \in \mathcal{S}^*$.*

PROOF.   Define an auxiliary random variable $Z$ such that $\mathbb{P}\{Z > u\} = e^{-H(t(u),u)}$. First, we show that $Z \in \mathcal{S}^*$. Since $\mathcal{S}^*$ is closed under tail equivalence, this implies that $W_\mu(T)$ and $M_\mu(T)$ are in $\mathcal{S}^*$ as well. According to Corollary 3.9 in [15], it suffices to show that the hazard rate $q_Z(u)$ of $Z$ is regularly varying of index $\nu \in (-1, 0)$. In fact, we will prove that $q_Z(u) = \frac{d}{du} H(t(u), u)$ is regularly varying of index $\beta - 1$. From the expression for $H(t(u), u)$ we obtain

$$q_Z(u) = t'(u)Q'(t(u)) + \frac{\mu^2}{2} t'(u) - \mu + \frac{u}{t(u)} - \frac{1}{2}\frac{u^2}{t(u)^2}t'(u),$$

where $t'(u)$ satisfies

$$Q''(t(u)t'(u)) = \frac{u}{t(u)^2} - \frac{u^2}{t(u)^3}t'(u).$$

From this equation, it can be shown, using T1, that there exists a constant $\kappa$ such that

$$t'(u) = \frac{1}{\mu} + (\kappa + o(1))L(u)u^{\beta-1}.$$

Using a similar calculation (see also [14]), one can show that

$$t(u) = \frac{u}{\mu} - \beta\mu^{1-\beta}u^{\beta}L(u)(1 + o(1)).$$

Combining the above equations, one obtains, after a tedious but straightforward computation, that $q_Z(u)$ is, indeed, regularly varying of index $\beta - 1$.

Thus, we conclude that $Z \in \mathcal{S}^*$. By Proposition 5.1, $W_\mu(T)$ and $Z$ are tail-equivalent. Since $\mathcal{S}^*$ is closed under tail equivalence, it follows that $W_\mu(T) \in \mathcal{S}^*$. □

A second question in the next section concerns the tail behavior of both $M_\mu(T^r)$ and $M_\mu(T)^r$, with $T$ an On-period. Clearly, these random variables have different distributions in general. However, the next proposition shows that they are tail-equivalent if $T$ is long-tailed:

PROPOSITION 3.3. *If $T \in \mathcal{L}$, then*

$$\mathbb{P}\{M_\mu(T^r) > u\} \sim \frac{\mathbb{E}\{M_\mu(T)\}}{\mu\mathbb{E}\{T\}}\mathbb{P}\{M_\mu(T)^r > u\}.$$

PROOF. Write

$$I := \mathbb{P}\{M_\mu(T^r) > u\}$$
$$= \mathbb{P}\{T^r > \tau(u)\} = \frac{1}{\mathbb{E}\{T\}}\mathbb{E}\{(T - \tau(u))\mathbb{1}_{\{T > \tau(u)\}}\}$$

and

$$II := \mathbb{P}\{M_\mu(T)^r > u\}$$
$$= \frac{1}{\mathbb{E}\{M_\mu(T)\}}\mathbb{E}\{(M(T) - u)\mathbb{1}_{\{T > \tau(u)\}}\}.$$

Now consider the following linear combination of the expectations in $I$ and $II$:

$$P(u) := \mu\mathbb{E}\{(T - \tau(u))\mathbb{1}_{\{T > \tau(u)\}}\} - \mathbb{E}\{(M_\mu(T) - u)\mathbb{1}_{\{T > \tau(u)\}}\}$$
$$= \mathbb{E}\{[\mu(T - \tau(u)) - (M_\mu(T) - u)]\mathbb{1}_{\{T > \tau(u)\}}\}.$$

Conditioning on $T$ and $\tau(u)$, we obtain

$$P(u) = \int_{y=0}^{\infty} d\mathbb{P}\{T \le y\}\int_{z=0}^{y} d\mathbb{P}\{\tau(u) \le z\}\mathbb{E}\{\mu(y - z) - (M_\mu(y) - u)|\tau(u) = z\}.$$

To simplify this expression, we investigate the distribution of $M_\mu(y) - u|\tau(u) = z$, with $0 < z \le y$. Note that

$$(3.5) \qquad M_\mu(y) = \sup\left\{M_\mu(z), W_\mu(z) + \sup_{z \le s \le y}[W_\mu(s) - W_\mu(z)]\right\}.$$

Since the condition $\tau(u) = z$ is a.s. equivalent to $W_\mu(z) = M_\mu(z) = u$, we get from (3.5) that, if $\tau(u) = z$,

$$(3.6) \qquad M_\mu(y) - u = \sup_{z \leq s \leq y} [W_\mu(s) - W_\mu(z)].$$

Using the strong Markov propery and the fact that $W_\mu(t), t \geq 0$, has stationary increments, we can conclude from (3.6) that $M_\mu(y) - u | \tau(u) = z$ is distributed as $M_\mu(y - z)$. Thus, we have

$$P(u) = \int_{y=0}^\infty d\mathbb{P}\{T \leq y\} \int_{z=0}^y d\mathbb{P}\{\tau(u) \leq z\} \mathbb{E}\{\mu(y - z) - M_\mu(y - z)\}.$$

Since $\mathbb{E}\{M_\mu(t)\} \geq \mu t$, we obtain

$$(3.7) \qquad \begin{aligned} |P(u)| &= |\mathbb{E}\{[\mu(T - \tau(u)) - (M_\mu(T) - u)]\mathbb{1}_{\{T > \tau(u)\}}\}| \\ &= \int_{y=0}^\infty d\mathbb{P}\{T \leq y\} \int_{z=0}^y d\mathbb{P}\{\tau(u) \leq z\} \mathbb{E}\{M_\mu(y - z) - \mu(y - z)\}. \end{aligned}$$

Using the definition of $M_\mu(t)$ and the self-similarity of Brownian motion, we have

$$\mathbb{E}\{(M_\mu(t) - \mu t)\} \leq \sqrt{t}\, \mathbb{E}\left\{\sup_{0 < s < 1} W(s)\right\} := \bar{W}\sqrt{t}.$$

Inserting this in (3.7) and integrating over $y$ and $z$, we obtain

$$(3.8) \qquad |P(u)| \leq \bar{W}\mathbb{E}\{\sqrt{T - \tau(u)}\, \mathbb{1}_{\{T > \tau(u)\}}\}.$$

Noting that $\sqrt{t} \leq 1/\varepsilon + \varepsilon t$ for any $\varepsilon > 0$ and any $t \geq 0$, we deduce from (3.8)

$$|P(u)| \leq \frac{1}{\varepsilon}\bar{W}\mathbb{P}\{T > \tau(u)\} + \varepsilon\bar{W}\mathbb{E}\{(T - \tau(u))\mathbb{1}_{\{T > \tau(u)\}}\}.$$

Since $T$ is long-tailed, $T^r$ is long-tailed as well. Thus, we have, as $u \to \infty$, because $\tau(u) \to \infty$ a.s., $\mathbb{P}\{T > \tau(u)\}/\mathbb{P}\{T^r > \tau(u)\} \to 0$. We conclude that for any $\varepsilon > 0$,

$$\limsup_{u \to \infty} \frac{|P(u)|}{\mathbb{P}\{T^r > \tau(u)\}} \leq \varepsilon\bar{W}\mathbb{E}\{T\},$$

that is, $|P(u)| = o(\mathbb{P}\{T^r > \tau(u)\})$. Using the definitions of $I$, $II$ and $P(u)$, we conclude that

$$\mathbb{P}\{M_\mu(T^r) > u\} \sim \frac{\mathbb{E}\{M_\mu(T)\}}{\mu\mathbb{E}\{T\}}\mathbb{P}\{M_\mu(T)^r > u\},$$

which completes the proof.  □

3.2. *Workload asymptotics.* In this section, we apply the results of the previous section to obtain the tail asymptotics of the workload $V_{X+Y}^c$. As mentioned earlier, we will follow the framework of Asmussen, Schmidli and Schmidt [4]. Recall that the increment process associated with $S(t) = X(t) + Y(t) - ct$ [with $X(t) \equiv W(t)$] is regenerative w.r.t. the delayed renewal process $\{Z_n, n \geq 0\}$ defined in Section 2. Thus, we consider the embedded process

$$S_n := X(Z_n) + Y(Z_n) - cZ_n =: U_0 + U_1 + \cdots + U_n.$$

Note that $S_n - S_0$, $n \geq 1$, is a random walk. Furthermore, define

$$M_0 := \sup_{0 < t < Z_0} S(t),$$

$$M_n = \sup_{Z_{n-1} < t < Z_n} [S(t) - S_{n-1}].$$

In order to obtain the asymptotics of $V_{X+Y}^c$, we apply the results of Section 3.2 of [4]. To check the assumptions stated there, we need the asymptotic behavior of the random variables $U_0$, $U_1$, $M_0$ and $M_1$. This is covered by the following lemma.

LEMMA 3.2. (i) *If* $T_{\mathrm{on}}$ *satisfies* T1, *then* $U_0, M_0 \in \mathcal{S}$, *and*

$$\mathbb{P}\{U_0 > u\} \sim \mathbb{P}\{M_0 > u\} \sim p\mathbb{P}\{W_{r-c}(T_{\mathrm{on}}^r) > u\}.$$

(ii) *If* $T_{\mathrm{on}}$ *satisfies* T1, *then* $U_1, M_1 \in \mathcal{S}^*$, *and*

$$\mathbb{P}\{U_1 > u\} \sim \mathbb{P}\{M_1 > u\} \sim \mathbb{P}\{W_{r-c}(T_{\mathrm{on}}) > u\}.$$

PROOF. We only prove the statement for $U_0$ and $M_0$ (the proof for $U_1$ and $M_1$ is similar, but easier). Recall the construction of the On-Off process given in Section 2. With a slight abuse of notation, we can write

$$U_0 \stackrel{d}{=} IW_{r-c}(T_{\mathrm{on}}^r) + (1 - I)(W_{-c}(T_{\mathrm{off}}^r) + W_{r-c}(T_{\mathrm{on}}^r)).$$

In this expression, all components are independent. Because $W_{-c}(T_{\mathrm{off}}^r) \leq \sup_{t>0} W_{-c}(t)$, this random variable is light-tailed. Second, since $T_{\mathrm{on}} \in \mathcal{L}$, we have $\mathbb{P}\{T_{\mathrm{on}} > u\} = o(\mathbb{P}\{T_{\mathrm{on}}^r > u\})$. This implies, using Lemma 3.1,

$$\mathbb{P}\{W_{r-c}(T_{\mathrm{on}}) > u\} \sim \mathbb{P}\{M_{r-c}(T_{\mathrm{on}}) > u\}$$
$$= \mathbb{P}\{T_{\mathrm{on}} > \tau(u)\}$$
$$= o(\mathbb{P}\{T_{\mathrm{on}}^r > \tau(u)\})$$
$$= o(\mathbb{P}\{M_{r-c}(T_{\mathrm{on}}^r) > u\})$$
$$= o(\mathbb{P}\{W_{r-c}(T_{\mathrm{on}}^r) > u\}).$$

By a straightforward application of Propositions 3.1–3.3, it follows that $M_{r-c}(T_{on}^r)$ is subexponential if $T_{on}$ satisfies T1. Thus, using standard properties of subexponential distributions, we conclude that

$$\mathbb{P}\{U_0 > u\} = p\mathbb{P}\{W_{r-c}(T_{on}^r) > u\} + (1-p)\mathbb{P}\{W_{-c}(T_{off}^r) + W_{r-c}(T_{on}) > u\}$$
$$\sim p\mathbb{P}\{W_{r-c}(T_{on}^r) > u\}.$$

Moving now to the tail behavior of $M_0$, note that (with a slight abuse of notation)

$$M_0 \le IM_{r-c}(T_{on}^r) + (1-I)\left(\sup_{t>0} B_{-c}(t) + M_{r-c}(T_{on})\right).$$

Hence, using a similar argument as above, we obtain

$$\mathbb{P}\{M_0 > u\} \overset{<}{\sim} p\mathbb{P}\{M_{r-c}(T_{on}^r) > u\} \sim p\mathbb{P}\{W_{r-c}(T_{on}^r) > u\}.$$

The asymptotic lower bound is trivial, since $M_0 \ge U_0$.  □

PROOF OF THEOREM 3.1.   Lemma 3.2 allows us to apply Corollary 3.2(ii) of [4], which yields

$$(3.9) \qquad \mathbb{P}\{V_{X+Y}^c > u\} \sim \mathbb{P}\{U_0 > u\} + \mathbb{P}\left\{\sup_{n\ge 1} S_n - S_0 > u\right\}.$$

The first term is covered by Lemma 3.2. To deal with the second term, note that in view of Lemma 3.2(ii), $[U_1^+]^r \in \mathcal{S}$. Using Veraverbeke's theorem [22], we then obtain

$$\mathbb{P}\left\{\sup_{n\ge 1} S_n - S_0 > u\right\}$$
$$\sim \frac{1}{-[\mathbb{E}\{W_{r-c}(T_{on})\} + \mathbb{E}\{W_{-c}(T_{off})\}]}$$
$$\qquad \times \int_u^\infty \mathbb{P}\{W_{r-c}(T_{on}) + W_{-c}(T_{off}) > v\}\,dv$$
$$\sim \frac{1}{-[\mathbb{E}\{W_{r-c}(T_{on})\} + \mathbb{E}\{W_{-c}(T_{off})\}]} \int_u^\infty \mathbb{P}\{M_{r-c}(T_{on}) > v\}\,dv$$
$$= \frac{1}{-[(r-c)\mathbb{E}\{T_{on}\} - c\mathbb{E}\{T_{off}\}]}\mathbb{E}\{M_{r-c}(T_{on})\}\mathbb{P}\{M_{r-c}^r(T_{on}) > u\}.$$

Finally, noting that

$$p = \frac{\mathbb{E}\{T_{on}\}}{\mathbb{E}\{T_{on}\} + \mathbb{E}\{T_{off}\}}, \qquad \rho = rp,$$

we conclude that

$$\mathbb{P}\Big\{\sup_{n\geq 1} S_n - S_0 > u\Big\} \sim \frac{p}{c-\rho}\frac{\mathbb{E}\{M_{r-c}(T_{\mathrm{on}})\}}{\mathbb{E}\{T_{\mathrm{on}}\}}\mathbb{P}\{M_{r-c}^r(T_{\mathrm{on}}) > u\}.$$

Thus, by (3.9) and Proposition 3.3, we obtain

$$\mathbb{P}\{V_{X+Y}^c > u\} \sim p\frac{r-\rho}{c-\rho}\mathbb{P}\{M_{r-c}(T_{\mathrm{on}}^r) > u\}.$$

Applying Lemma 3.1 completes the proof. $\square$

**4. Oscillatory behavior.** In this section we assume that $X(t)$ is Gaussian satisfying C1 and C2 and that $Y(t)$ is an integrated On-Off process with peak rate $r$. The main difference from the previous section is that we now assume that $r = c$. Under this critical condition, the process $S(t) = X(t) + Y(t) - ct$ will oscillate during the On-periods of $Y$. The next theorem gives the main result of this section.

THEOREM 4.1. *If $X(t)$ has stationary increments and satisfies conditions* C1 *and* C2, *and $Y(t)$ is an integrated On-Off process with $T_{\mathrm{on}}$ regularly varying of index $-\nu < -1$ such that $\mathbb{P}\{T_{\mathrm{on}}^r > u\} = L(x)u^{1-\nu}$ with $L(\cdot)$ slowly varying, then*

$$\mathbb{P}\{V_{X+Y}^r > u\} \sim p\mathbb{E}\{\bar{B}_H^{H(\nu-1)}\}\mathbb{P}\{\sigma_X(T_{\mathrm{on}}^r) > u\},$$

*with $H = \alpha/2$ and $\bar{B}_H = \sup_{0\leq s\leq 1} B_H(s)$. In particular, $V_{X+Y}^r$ is regularly varying of index $(1-\nu)/H$.*

Theorem 4.1 shows that the heaviness of the tail of $V_{X+Y}^r$ is a combined effect of the heaviness of $T_{\mathrm{on}}^r$ and the degree of dependence in $X$. Informally, a large value of $V_{X+Y}^c$ is most likely caused by a single long On-period which started at time 0. During this long On-period, the net input process has zero drift. This implies that the net input process at time $t$ is $O(\sigma_X(t))$. Hence, in order to reach level $u$, the length of the On-period has to be $O(\sigma_X^{-1}(u))$.

In the proof of Theorem 4.1, we make these heuristics precise. Before we give a proof, we first present some auxiliary results. The following lemma is taken from [11].

LEMMA 4.1. *Under the conditions of Theorem 4.1, we have*

$$\mathbb{P}\Big\{\sup_{0\leq t\leq T_{\mathrm{on}}^r} X(t) > u\Big\} \sim \mathbb{E}\{\bar{B}_H^{H(\nu-1)}\}\mathbb{P}\{\sigma_X(T_{\mathrm{on}}^r) > u\}.$$

The main idea of the proof of Theorem 4.1 is to separate the processes $X$ and $Y$ by adding and subtracting nonlinear perturbations. To handle such perturbations, we need a further auxiliary lemma, which is Proposition 1 of [12].

LEMMA 4.2. *Let $X(t)$ be a centered Gaussian process satisfying conditions* C1 *and* C2. *If $\eta > \alpha/2$, then*

$$\log \mathbb{P}\Big\{ \sup_{t \geq 0} [X(t) - dt^\eta] > u \Big\} \sim -\frac{1}{2} d^{\alpha/\eta} \Big( \frac{\alpha}{2\eta - \alpha} \Big)^{-\alpha/\eta} \Big( \frac{2\eta}{2\eta - \alpha} \Big)^2 \frac{u^2}{\sigma_X^2(u^{1/\eta})}.$$

PROOF OF THEOREM 4.1.   The lower bound is trivial, in view of Lemma 4.1 and the construction of the process $Y(t)$ given in Section 2.

For the upper bound, write for some $\gamma \in (0, 1)$,

$$\mathbb{P}\Big\{ \sup_{t \geq 0} S(t) > u \Big\} \leq \mathbb{P}\Big\{ \sup_{t \leq Z_0} S(t) > (1 - \gamma)u \Big\} + \mathbb{P}\Big\{ \sup_{t > Z_0} S(t) - S(Z_0) > \gamma u \Big\}.$$

We need to show that the second term can be asymptotically neglected. Using sample-path arguments, we have

$$\mathbb{P}\Big\{ \sup_{t > Z_0} S(t) - S(Z_0) > \gamma u \Big\}$$

$$= \mathbb{P}\Big\{ \sup_{t > Z_0} [X(t) - X(Z_0) + Y(t) - Y(Z_0) - r(t - Z_0)] > \gamma u \Big\}$$

$$\leq \mathbb{P}\Big\{ \sup_{t > Z_0} [Y(t) - Y(Z_0) - r(t - Z_0) + d(t - Z_0)^\eta] > \gamma u / 2 \Big\}$$

$$\quad + \mathbb{P}\Big\{ \sup_{t > Z_0} [X(t) - X(Z_0) - d(t - Z_0)^\eta] > \gamma u / 2 \Big\}$$

$$= I + II,$$

where we take $1 > \eta > \alpha/2$ and $d$ small.

We first deal with term $I$. Observe that

$$d(Z_n - Z_0)^\eta \leq d \sum_{i=1}^n (Z_i - Z_{i-1})^\eta,$$

from which it follows that

$$I \leq \mathbb{P}\Big\{ \sup_{n \geq 1} S_n > \gamma u / 2 \Big\},$$

where $S_n$ is a random walk with generic step size $U = dT_{\mathrm{on}}^\eta + dT_{\mathrm{off}}^\eta - rT_{\mathrm{off}}$.

We can choose $d$ small enough such that $U$ has negative mean. Noting that $dT_{\mathrm{off}}^\eta - rT_{\mathrm{off}}$ is bounded from above, we conclude that the right tail of $U$ is regularly varying. This allows us to apply Veraverbeke's theorem [22], yielding

$$I \leq \mathbb{P}\Big\{ \sup_{n \geq 1} S_n > \gamma u / 2 \Big\} \sim \frac{1}{-\mathbb{E}\{U\}} \int_{\gamma u / 2}^\infty \mathbb{P}\{U > y\} \, dy,$$

which is regularly varying of index $1 - \nu\eta$. We can choose $\eta$ such that $1 - \nu\eta > (1 - \nu)H$ (i.e., $\eta < H + \frac{1-H}{\nu}$).

We now turn to term $II$. This term is somewhat easier: since $X(t)$ has stationary increments, we have

$$II = \mathbb{P}\Big\{\sup_{t \geq 0}[X(t) - dt^\eta] > x\Big\}.$$

This probability is decaying faster than any polynomial, in view of Lemma 4.2. Thus, we can conclude that, for any $\gamma > 0$,

(4.1)
$$\mathbb{P}\Big\{\sup_{t \geq 0} S(t) > u\Big\} \underset{\sim}{<} \mathbb{P}\Big\{\sup_{t \leq Z_0} S(t) > (1 - \gamma)u\Big\}.$$

We finally evaluate the probability on the right-hand side by conditioning upon the state of the On-Off process $J$ at time 0

$$\mathbb{P}\Big\{\sup_{t \leq Z_0} S(t) > (1 - \gamma)u\Big\}$$

$$= p\mathbb{P}\Big\{\sup_{t \leq T_{\mathrm{on}}^r} X(t) > (1 - \gamma)u\Big\}$$

$$+ (1 - p)\mathbb{P}\Big\{\sup_{t \leq T_{\mathrm{off}}^r + T_{\mathrm{on}}}[Y(t) + X(t) - rt] > (1 - \gamma)u\Big\}.$$

Using similar methods as above, it is straightforward to show that the second term is regularly varying of index $-\nu/H$. From the proof of the lower bound, we already know that the first term is regularly varying of index $(1 - \nu)/H$. Hence, we conclude from (4.1) and Lemma 4.1,

$$\limsup_{u \to \infty} \frac{\mathbb{P}\{\sup_{t \geq 0} S(t) > u\}}{p\mathbb{E}\{\bar{B}_H^{H(\nu-1)}\}\mathbb{P}\{T_{\mathrm{on}}^r > \sigma_X(u)\}} \leq (1 - \gamma)^{-(1-\nu)/H}$$

for all $\gamma > 0$. $\quad\square$

We conclude this section by noting that the pre-factor $\mathbb{E}\{\bar{B}_H^{H(\nu-1)}\}$ can be computed explicitly when $H = 1/2$, or, equivalently, when $\alpha = 1$:

COROLLARY 4.1. *In addition to the assumptions of Theorem 4.1, assume that $X(t)$ satisfies conditions C1 and C2, with $\alpha = 1$. Then*

$$\mathbb{P}\{V_{X+Y}^r > u\} \sim p\frac{1}{\sqrt{\pi}}2^{1+\nu}\Gamma\Big(\nu + \frac{1}{2}\Big)\mathbb{P}\{\sigma_X(T_{\mathrm{on}}^r) > u\}.$$

PROOF. The result follows in a straightforward manner from Theorem 4.1, combined with Proposition 2.1 in [11]. $\quad\square$

**5. Large deviations: reduced-peak equivalence.** In this section we consider the case that $X$ is Gaussian and $Y$ is an integrated On-Off process with peak rate $r < c$. We assume that the tail of $V_Y^d$, $\rho < d < r$, is heavier than that of $V_X^c$.

Under these conditions, it is clear that a reduced-load equivalence cannot hold. Informally, one can observe that $X$ cannot be replaced by its mean $(0)$, since $V_Y^c \equiv 0$, nor can $Y$ be replaced by its mean, since it has heavier tails than $X$. In fact, the next theorem shows that *both* $X$ and $Y$ need to show atypical behavior in order for the process $S(t) = X(t) + Y(t) - ct$ to reach a large value.

THEOREM 5.1. *Suppose that the process $X(t)$ has stationary increments and satisfies conditions* C1 *and* C2. *Furthermore, let $Y(t)$ be an integrated On-Off process, with $T_{\mathrm{on}}^r$ regularly varying, and $r < c$. Then*

$$\mathbb{P}\{V_{X+Y}^c > u\} \sim p\mathbb{P}\{V_X^{c-r} > u\}\mathbb{P}\Big\{T_{\mathrm{on}}^r > \frac{1}{c-r}\frac{\alpha}{2-\alpha}u\Big\}.$$

Theorem 5.1 may be combined with results in [10] or [17] to obtain an explicit expression for the asymptotic behavior of $\mathbb{P}\{V_{X+Y}^c > u\}$. A similar "reduced-peak" equivalence result has been proved in Theorem 3.1 of [7] for the case where $X$ is not a Gaussian process, but a general light-tailed process satisfying a large-deviations principle.

PROOF OF THEOREM 5.1. Let $t_u = \frac{1}{d}\frac{\alpha}{2-\alpha}u$. To let Theorem 3.1 of [7] carry over to the setting of the present paper, it is sufficient to prove analogues of Propositions 2.1 and 2.2 of [7]. In our setting, these propositions state that, for every $d > 0$, as $u \to \infty$,

$$(5.1) \qquad \frac{u^\beta \mathbb{P}\{V_X^{d+\delta} > u\}}{\mathbb{P}\{V_X^d > u\}} \to 0, \qquad \beta < \infty,$$

$$(5.2) \qquad \frac{u^\beta \mathbb{P}\{V_X^d([0,(1-\varepsilon)t_u]) > u\}}{\mathbb{P}\{V_X^d > u\}} \to 0, \qquad \varepsilon > 0, \beta < \infty,$$

$$(5.3) \qquad \frac{\mathbb{P}\{V_X^d([0,(1+\varepsilon)t_u]) > u\}}{\mathbb{P}\{V_X^d > u\}} \to 1, \qquad \varepsilon > 0.$$

Thus, it is sufficient to prove (5.1)–(5.3). To prove (5.1), note that Lemma 4.2 with $\eta = 1$ implies that, for some constant $C = C_{d,\delta} > 0$,

$$\frac{\mathbb{P}\{V_X^{d+\delta} > u\}}{\mathbb{P}\{V_X^d > u\}} = e^{-C(1+o(1))u^2/\sigma_X^2(u)},$$

which implies (5.1).

To prove (5.3), we define $\tau(u) := \inf\{t : X(t) - dt = u\}$ and note that Theorem 1 in [12] implies that $\tau(u)/t_u \to 1$ in $\mathbb{P}\{\cdot|\tau(u) < \infty\}$-probability.

It remains to prove (5.2). For this, define $X_u(t) = X(t)/(u + dt)$. Using the Borell inequality ([1], page 43), we obtain, for all $u > 0$,

$$
\begin{aligned}
&\mathbb{P}\{V_X^d([0, (1-\varepsilon)t_u]) > u\} \\
&= \mathbb{P}\left\{\sup_{t \in [0, (1-\varepsilon)t_u]} X_u(t) > 1\right\} \\
&\leq 2\exp\left(-\left(1 - \mathbb{E}\left\{\sup_{t \geq 0} X_u(t)\right\}\right)^2 \min_{t \leq (1-\varepsilon)t_u} \frac{(u+dt)^2}{2\sigma_X^2(t)}\right).
\end{aligned}
$$

Since $\lim_{u \to \infty} \mathbb{E}\{\sup_{t \geq 0} X_u(t)\} = 0$, by Lemma 2.2 in [9], we have

$$
\log(\mathbb{P}\{V_X^d([0, (1-\varepsilon)t_u]) > u\}) \underset{\sim}{<} -\min_{t \leq (1-\varepsilon)t_u} \frac{(u+dt)^2}{2\sigma_X^2(t)}.
$$

Using the uniform-convergence theorem for regularly varying functions, we obtain

$$
\begin{aligned}
\lim_{u \to \infty} \frac{\sigma_X(u)^2}{u^2} \min_{t \leq (1-\varepsilon)t_u} \frac{(u+dt)^2}{2\sigma_X^2(t)} &= \lim_{u \to \infty} \min_{s \leq (1-\varepsilon)\alpha/(d(2-\alpha))} \frac{(1+ds)^2}{2\sigma_X^2(su)/\sigma_X^2(u)} \\
&= \min_{s \leq (1-\varepsilon)\alpha/(d(2-\alpha))} \frac{(1+ds)^2}{2s^\alpha} \\
&= d^\alpha \frac{(1+(1-\varepsilon)\alpha/(2-\alpha))^2}{2(1-\varepsilon)^\alpha} \left(\frac{\alpha}{2-\alpha}\right)^\alpha \\
&> 2d^\alpha \frac{(2-\alpha)^{\alpha-2}}{\alpha^2},
\end{aligned}
$$

where the last inequality is valid for all $\varepsilon > 0$.

Now note that, in view of Lemma 4.2,

$$
\log(\mathbb{P}\{V_X^d > u\}) \sim -2d^\alpha \frac{(2-\alpha)^{\alpha-2}}{\alpha^2} u^2/\sigma_X^2(u).
$$

Putting everything together, we conclude that, for every $\varepsilon > 0$, there exists a constant $K_\varepsilon$ such that

$$
\frac{\mathbb{P}\{V_X^d([0, (1-\varepsilon)t_u]) > u\}}{\mathbb{P}\{V_X^d > u\}} \leq e^{-K_\varepsilon(1+o(1))u^2/\sigma_X^2(u)}.
$$

This implies (5.2). $\square$

**Acknowledgments.** The authors would like to thank Ton Dieker for showing us his recent work [12] and an anonymous referee for providing useful comments.

B. ZWART
DEPARTMENT OF MATHEMATICS
  AND COMPUTER SCIENCE
EINDHOVEN UNIVERSITY OF TECHNOLOGY
P.O. BOX 513
5600 MB EINDHOVEN
THE NETHERLANDS
E-MAIL: zwart@win.tue.nl

S. BORST
CWI
P.O. BOX 94079
1090 GB AMSTERDAM
THE NETHERLANDS
E-MAIL: sem@cwi.nl

K. DĘBICKI
MATHEMATICAL INSTITUTE
UNIVERSITY OF WROCŁAW
PL. GRUNWALDZKI 2/4
50-384 WROCŁAW
POLAND
E-MAIL: debicki@math.uni.wroc.pl